\documentclass{article}

\usepackage{PRIMEarxiv}

\usepackage{abstract}
\usepackage[utf8]{inputenc} 
\usepackage[T1]{fontenc}   
\usepackage{hyperref}       
\usepackage{url}            
\usepackage{booktabs}       
\usepackage{amsfonts}       
\usepackage{amsmath}        
\usepackage{amsthm}         
\usepackage{nicefrac}       
\usepackage{microtype}      
\usepackage{lipsum}
\usepackage{fancyhdr}       
\usepackage{graphicx}       
\usepackage{subcaption}     
\graphicspath{{media/}}     

\newtheorem{theorem}{Theorem}[section]

\title{Prime Factorization Equation from a Tensor Network Perspective
 }

\author{
    Alejandro Mata Ali \\
    Instituto Tecnológico de Castilla y León, Burgos, Spain\\
    \texttt{alejandro.mata@itcl.es} \\
     \And
    Jorge Martínez Martín \\
    Instituto Tecnológico de Castilla y León, Burgos, Spain\\
    \texttt{jorge.martinez@itcl.es} \\
     \And
    Sergio Muñiz Subiñas \\
    Instituto Tecnológico de Castilla y León, Burgos, Spain\\
    \texttt{sergio.muniz@itcl.es} \\
     \And
    Miguel Franco Hernando \\
    Instituto Tecnológico de Castilla y León, Burgos, Spain\\
    \texttt{miguel.franco@itcl.es} \\
     \And
    Javier Sedano \\
    Instituto Tecnológico de Castilla y León, Burgos, Spain\\
    \texttt{javier.sedano@itcl.es} \\
     \And
    Ángel Miguel García-Vico \\
    Andalusian Research Institute in Data Science and Computational Intelligence (DaSCI),\\
    University of Jaén, 23071 Jaén, Spain\\
    \texttt{agvico@ujaen.es} \\
  }

\begin{document}
\maketitle

\begin{abstract}
This paper presents an exact and explicit tensor-network equation for the search of nontrivial divisors of a composite integer, together with an algorithm for its computation. The proposed method is based on the MeLoCoToN approach, which addresses combinatorial optimization problems through classical tensor networks. The presented tensor network tensorizes a binary multiplication circuit and projects its output onto the target integer to be factorized. Additionally, in order to make the algorithm more efficient, the number and dimension of the tensors and their contraction scheme are optimized, including a reduced auxiliary register that still preserves at least one valid factorization orientation. Finally, a series of tests on the algorithm are conducted, contracting the tensor network both exactly and approximately using tensor train compression, and evaluating its performance.
\end{abstract}

\tableofcontents

\section{Introduction}
Prime factorization is the process of expressing a positive composite number $N$ as the product of its prime factors; that is, finding all $p_i$ that satisfy the equation $N = \prod_i p_i$. Despite its apparent simplicity, the problem is very difficult to solve. It has efficiently verifiable certificates, and standard decision formulations lie in NP and coNP, but no classical polynomial-time factoring algorithm is known and the problem is not known to be NP-complete. The problem has been studied for centuries due to its academic significance~\cite{Prime_Fact_Survey} and its foundational role in areas such as number theory or cryptography. This is of special interest for semiprime numbers, which are products of two prime numbers, so $N=pq$. This has led to the development of various prime factorization algorithms, including memcomputing techniques~\cite{MEMCPU_Prime} and genetic algorithms~\cite{Genetic_Prime}. Consequently, the security of the widely used RSA~\cite{RSA} cryptographic algorithm depends on the difficulty of this problem, as deriving the private key from the public key requires the factorization of large numbers. This is an important topic in the cybersecurity field, and even more so in conflict contexts. Large RSA challenge numbers such as RSA-240 have been factored using highly optimized variants of the general number field sieve~\cite{RSA_240}.

In this problem, quantum computing offers an advantage with Shor's algorithm for prime factorization~\cite{Shor,Shor_2n_3,Shor_Toffoli}, which has an exponential acceleration regarding the best known classical algorithms. In addition, there are other quantum algorithms, based on QUBO formulation~\cite{HUBO_Prime_Fact} or variational circuits~\cite{VQE_Prime} that allow its resolution. Nevertheless, these algorithms are limited by the quality of current error-prone quantum hardware, which is why several alternative quantum-inspired technologies have been developed. One of the most interesting approaches is the use of tensor networks~\cite{TN}, consisting of the graphical representation of a tensor contraction equation. Unlike quantum computing algorithms, the tensor networks are computed in classical resources because it is a classical computation approach. This technology allows for efficient classical simulation of some low-entangled quantum systems by transforming them into tensor representations.  An example is the Tensor Train (TT) representation~\cite{TT}, which consists of a linear chain of tensors connected to each other by their bond indices. These representations are useful for computing tensor network contractions, since they allow approximating the state at each step without requiring excessive scaling. This approach has been explored to break up to RSA-100~\cite{RSA_TN}. Another possible approach is the MeLoCoToN formalism~\cite{melocoton} which allows obtaining the exact and explicit equation that solves inversion and optimization problems by using classical logical circuits and tensor networks. This approach has been applied to several optimization and constraint satisfaction problems, and as far as is known, never to the factorization case. A similar approach is presented in~\cite{constrain} for creating tensor networks that implement exactly arbitrary discrete linear constraints, more focused on tensor trains.

The main contributions of this work are as follows:
\begin{itemize}
    \item An exact and explicit tensor-network equation that encodes the search for a nontrivial divisor of a given odd semiprime number, with a direct extension to more general composite inputs.
    \item An optimization of this equation to reduce computational resources, particularly through the reduction of the auxiliary $q$-register and in the case of RSA-like instances.
    \item Two tensor network-based algorithms to compute a valid factorization result and an analysis of their computational complexity.
    \item An study of the computability of the correct result with compressions of tensor trains.
\end{itemize}

This work will be structured as follows. First, Sec.~\ref{sec: background} introduces a brief background on the semiprime factorization algorithms. Secondly, Sec.~\ref{sec: tensor networks} introduces the tensor network equation and algorithm. Then, Sec.~\ref{sec: performance} presents the experiments performed and their results. Finally, Sec.~\ref{sec: discussions} discusses the results.

\section{Background}\label{sec: background}
In the field of semiprime factorization, there exist several algorithms and studies in the literature, due to their relevance. This section provides a brief background of the available algorithms, focusing on the most relevant ones. This work focuses on odd numbers, because if $N$ is even then the nontrivial factor $2$ is obtained immediately and one may continue recursively on $N/2$.

The first method is the Fermat factorization method, which consists of looking for a value $a$, which satisfies
\begin{equation}
    N = a^2 -b^2 = (a+b)(a-b).
\end{equation}
In this way, if $N=pq$, then
\begin{equation}
    N=\left(\frac{p+q}{2}\right)^2 - \left(\frac{p-q}{2}\right)^2.
\end{equation}

The algorithm starts from an initial value $a=\left\lceil \sqrt{N}\right\rceil$, and then calculates $b_2=a^2-N$. This process of increasing $a$ by one unit is repeated until $b_2$ is a square. If $N$ is a prime, it requires $O(N)$ steps, but if both numbers are similar, near $\sqrt{N}$, it is a fast algorithm to reach factorization. If it is combined with the trial division~\cite{division_trial}, its complexity is reduced to $O(N^{1/3})$ steps, which remains exponential with respect to the number of bits of $N$.

The second algorithm is Dixon's factorization method~\cite{Dixon_Fact}. It is based on finding a congruence of squares modulo $N$ such as Fermat's method. A bound $B$ is imposed in order to define a base $P$ of all primes equal to or lower than $B$. Then a set of values of $z$ are selected, such that $z^2 \mod N$ is $B$-smooth. With these values and linear algebra, it reaches the congruent pair of values $a$ and $b$. Its complexity is $O\left(\exp\left(2\sqrt{2}\sqrt{\ln N \ln \ln N}\right)\right)$.

The third algorithm is the Quadratic Sieve~\cite{quadratic_sieve}, the second-fastest method known in practice, being the fastest for integers under 100 decimal digits. It is a modification of Dixon's method, but it makes use of quadratic polynomials. Its complexity is conjectured to be\\ \mbox{$O(\exp((1+O(1))\sqrt{\ln N \ln \ln N}))$}.

The last classical method is the General Number Field Sieve~\cite{GNFS}, which is the most efficient classical algorithm known for this task for numbers larger than $10^{100}$. Its complexity is $O\left(\exp \left( ((64/9)^{1/3}+O(1) ) (\ln N)^{1/3}(\ln \ln N)^{2/3} \right)\right)$.

In the quantum computing field, the most well-known algorithm is Shor's algorithm~\cite{Shor}, which reduces the problem to looking for the period of a modular function. Its complexity is $O((\ln N)^2(\ln\ln N)(\ln\ln\ln N))$, polynomial in the number of bits of $N$~\cite{Shor_complex}.

\section{Tensor Network for Factorization}\label{sec: tensor networks}
This section introduces the tensor network equation and the algorithm to compute it. The MeLoCoToN formalism~\cite{melocoton} is used to obtain the equation, and then design an optimized algorithm to compute it. The MeLoCoToN formalism has three basic steps. First, it needs the creation of the classical logical circuit that computes the function to invert. In Sec.~\ref{ssec: logical} the creation of the logical circuit for the multiplication relation is described, with the odd semiprime setting used as the main motivating case. The second step is to tensorize the logical circuit to obtain the tensor network that represents the logical space. The last step is to iterate on this tensor network, contracting it with other tensors to determine a consistent factor value and getting the explicit equation. The last two steps are performed in Sec.~\ref{ssec: tensorization}, optimizing the size of the tensors involved. Finally, to obtain the solution values, it is necessary to contract these tensor networks, which requires a contraction scheme. In Sec.~\ref{ssec: complexity} two contraction schemes are provided, together with an analysis of their computational complexity.

\subsection{Logical Circuit}\label{ssec: logical}
This subsection introduces the classical logical circuit that serves as the foundation for the tensor network. This logical circuit takes $p$ and $q$ and represents the relation $y=pq$, where $y$ is later projected onto the target integer $N$. For pedagogical reasons, the construction starts from the standard modular multiplication circuit and then adapts it to the exact factorization setting. In the final formulation, if $p$ is represented with up to $n-1$ bits and the reduced auxiliary register $q$ uses $m=\lceil n/2\rceil$ bits, the output register is chosen wide enough to represent the full product, for example with $W=n+m$ bits. The target integer $N$ is then padded with leading zeros on this register, so the final projection enforces exact integer multiplication rather than equality modulo a power of two.

To provide a clear understanding of the logical circuit, it is convenient to begin by explaining the logic behind an adder circuit of two numbers in Sec.~\ref{sssec: add} and a modular binary multiplier of two numbers in Sec.~\ref{sssec: multiplication}. These are known operations in binary computing, but the former will be adapted to this problem. In these cases, the functions $f_q(p) = p +q$ and $g_q(p)=pq$ are computed, respectively. These functions receive only one of the numbers and then operate it with the other number, which is a parameter of the function. To create functions of two numbers, $f(p,q)=p+q$ and $g(p,q)=pq$, these types of circuits need to be generalized, as in Sec.~\ref{sssec: double}, where the final exact multiplication relation is enforced.

\subsubsection{Binary addition}\label{sssec: add}
\begin{figure*}
    \centering
    \includegraphics[width=0.75\linewidth]{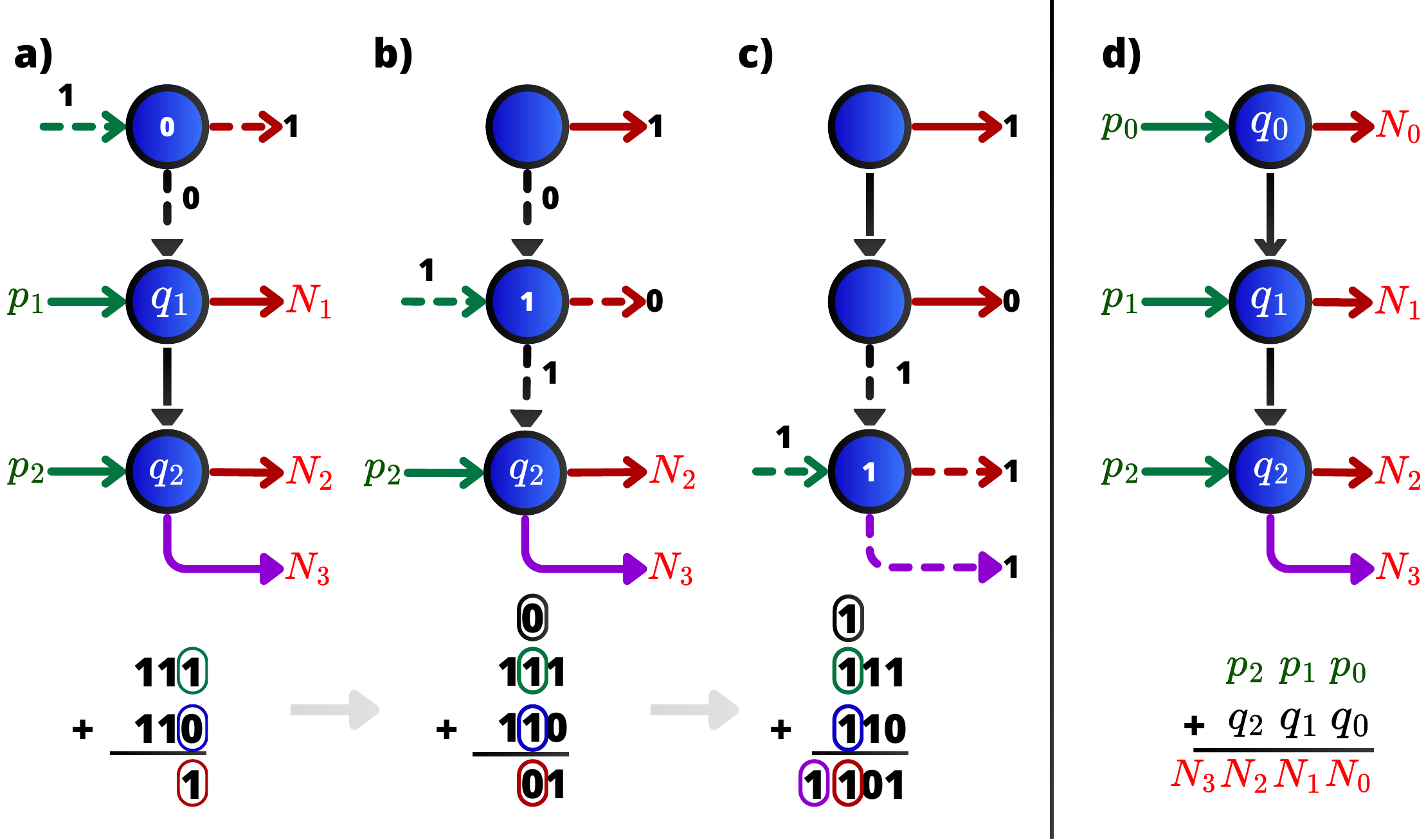}
    \caption{Logical circuit of the binary addition between two numbers. a), b) and c) represent each step of a binary addition between the numbers p = 111 and q = 110. d) Abstraction of the logical circuit of a binary addition of two numbers with 3 bits at most.}
    \label{fig: Binary Addition}
\end{figure*}

To illustrate the binary adder to add an input $p$ and a fixed $q$, consider a specific example where $p = p_2p_1p_0=111$ and $q=q_2q_1q_0=110$, shown in Fig.~\ref{fig: Binary Addition}. Starting with the least significant bit, a bitwise addition is performed, resulting in the carry $$c_0 = \left\lfloor\frac{p_0+q_0}{2}\right\rfloor=0$$ and the sum $$N_0 = p_0\oplus q_0=1$$ shown in Fig.~\ref{fig: Binary Addition} a. The next operator receives $c_0$, $p_1$ and $q_1$ and computes $$c_1 = \left\lfloor\frac{c_0+p_1+q_1}{2}\right\rfloor=1$$ and $$N_1 = c_0\oplus p_1\oplus q_1=0,$$ as shown in Fig.~\ref{fig: Binary Addition} b. The last operator receives $c_1$, $p_2$ and $q_2$ and computes
$$c_2 = \left\lfloor\frac{c_1+p_2+q_2}{2}\right\rfloor=1 =N_3$$
and
$$N_2 = c_1\oplus p_2\oplus q_2=1,$$ shown in Fig.~\ref{fig: Binary Addition} c. If final value $c_2$ is neglected, it performs the modular binary addition \mbox{$(p+q) \mod 2^n$}.

\begin{figure}
    \centering
    \includegraphics[width=0.75\linewidth]{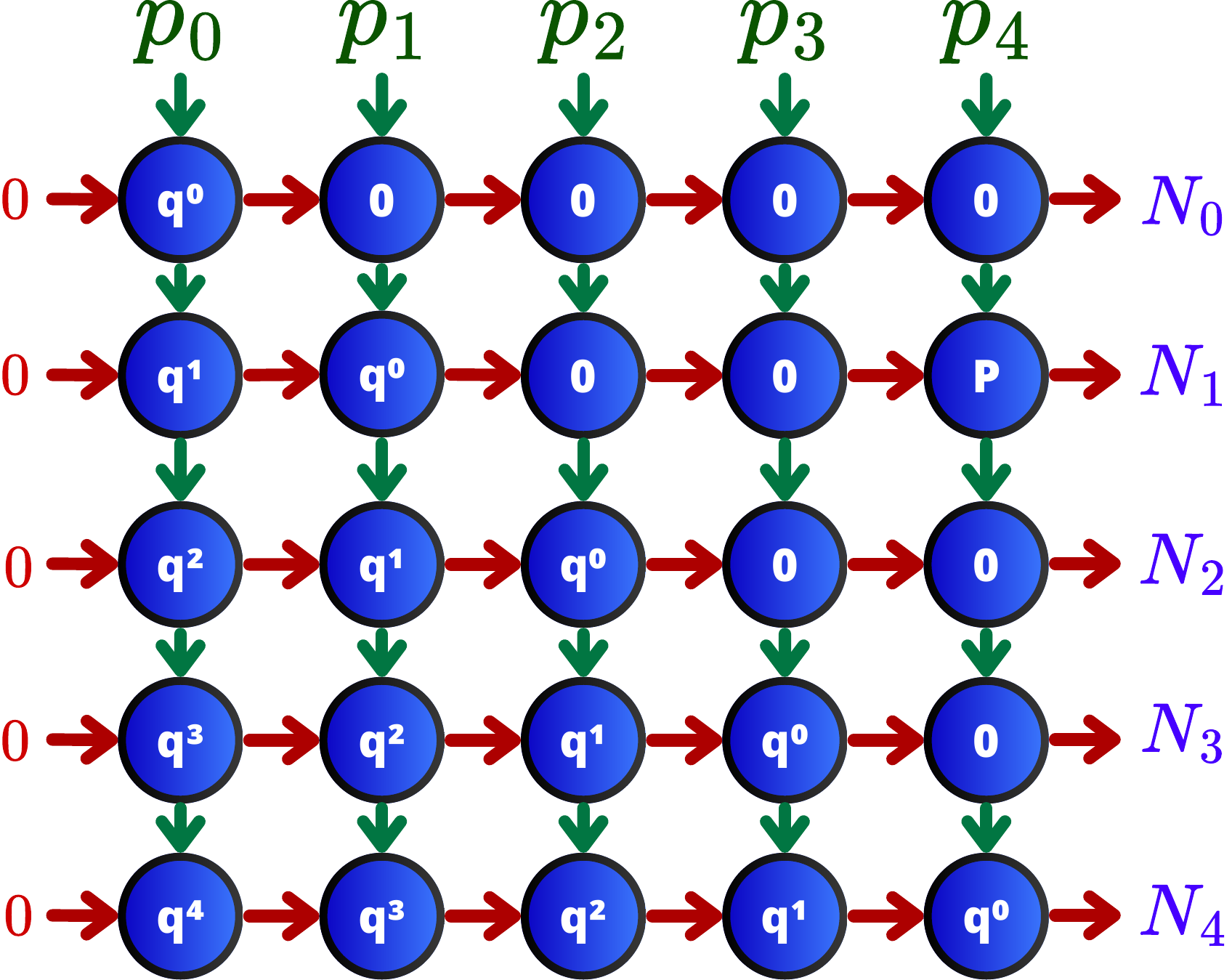}
    \caption{Logical circuit of the binary multiplication between two numbers $p$ and $q$.}
    \label{fig: Binary Multiplication}
\end{figure}

\subsubsection{Modular binary multiplication}\label{sssec: multiplication}
The logical circuit presented in Fig.~\ref{fig: Binary Multiplication} represents the modular binary multiplication of an input $p$ and a fixed $q$ based on Eq.~\ref{eq: multiplication}. This logical circuit can be interpreted as a series of controlled additions. Taking into account the following property~\cite{Shor_2n_3}
\begin{equation}\label{eq: multiplication}
    (qp)\text{ mod }2^n = (...((2^0 qp_0)\text{ mod }2^n + 2^1 qp_1)\text{ mod }2^n+... +2^{n-1} qp_{n-1})\text{ mod }2^n,
\end{equation}
    
a modular multiplication can be performed by an iterative modular summation. Starting from zero, a concatenated set of modular additions of Fig.~\ref{fig: Binary Addition} are performed, each one adding \mbox{$(2^k q)\mod 2^n$} only if $p_k=1$. To simplify the following computations, it must be considered that $(2^k q)\mod 2^n$ is the same binary string as $q$, but neglect the $k$ most significant bits and concatenate $k$ zeros to its end. An example is that \mbox{$(2\times5)\mod 2^3=2$} is in binary \mbox{$(2\times 101)\mod 2^3 = 010$}. Then, the same addition layer can be used repeatedly, removing the operators of the neglected bits.

\subsubsection{Binary multiplication of two numbers}\label{sssec: double}
\begin{figure*}
    \centering
    \includegraphics[width=\linewidth]{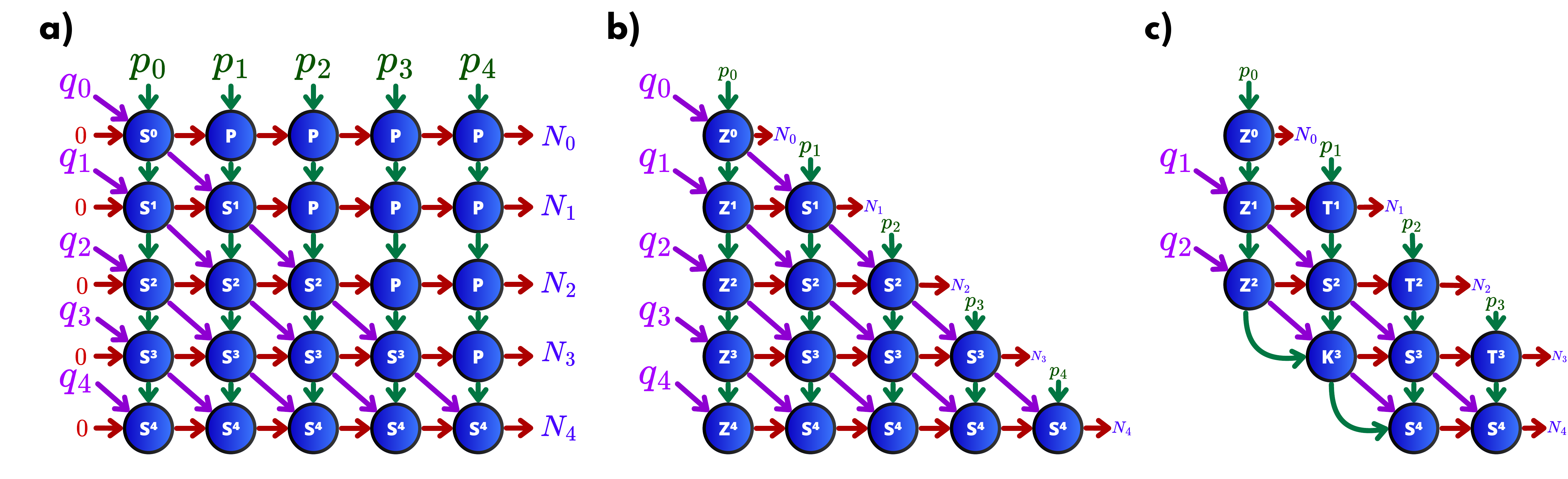}
    \caption{a) Logical circuit that performs the multiplication in binary terms of two numbers $p$ and $q$. In the exact formulation, the output register is understood as wide enough for the full product, with the high-order bits projected to zero when the target $N$ has only $n$ significant bits. b) Logic circuit that performs the same calculation, eliminating the unnecessary $P$ operators and $b_0=0$ input. c) Logic circuit in which the number $q$ has up to $\lceil n/2\rceil$ bits and the number $p$ has up to $n-1$ bits.}
    \label{fig: Logical Circuit}
\end{figure*}

In the previous section, the logical circuit that computes $(qp)\text{ mod }2^n$ for an input $p$ and a fixed $q$ has been shown. However, some modifications are necessary to achieve the goal of designing a circuit capable of taking as input both $p$ and $q$ and then enforcing exact multiplication. Fig.~\ref{fig: Logical Circuit} a illustrates the generalization so that each layer $k$ receives as horizontal input $b_{k}$ (Eq.~\ref{sum}) and vertical input $p_k$. 
\begin{equation}\label{sum}
    b_{k+1} = b_k + 2^{k}q p_k = \sum_{j=0}^k 2^j q p_j,\qquad b_0 = 0.
\end{equation}
To ensure that the application of each $q$ is consistent throughout all layers, it is sent as a diagonal signal, since the binary of $2x$ is the same number $x$ shifted a bit. For the factorization problem, the output register is understood as wide enough to represent the full product. Equivalently, one may regard the optimized circuit as projecting the bits of $N$ on its $n$ least significant output wires and projecting the remaining high-order output wires to zero. In this way, the network keeps only the assignments compatible with the exact integer equality $pq=N$.

As $q$ is shifted one bit per layer, the first $k$ operators of the $k$-th layer are operators that pass the control value $p_k$ and the corresponding bit of $b$ without modifying it. For this reason, these operators do not add any information and can be removed, leaving only the operators of the lower triangle (Fig. \ref{fig: Logical Circuit} b). Moreover, since the initial value $b_0$ is always $0$, the first column of operators can be directly applied to sum on the value $0$, so that they do not receive an input $b$.

Let $N$ be an $n$-bit composite integer and let $m=\lceil n/2\rceil$. Since $N$ is composite, it has a nontrivial divisor $d$ with $1<d\leq \sqrt{N}$. Because $\sqrt{N}<2^{n/2}\leq 2^m$, the divisor $d$ can be represented in the reduced $q$-register. Therefore, the orientation $q=d$ and $p=N/d$ is representable. This reduced $q$-register is an auxiliary device ensuring that at least one valid factorization orientation is representable. It should not be confused with an ordering constraint, and it does not force the exposed $p$-register to contain the smaller or the larger factor. In this way, half of the lower diagonals are removable since they only add up $0$, which implies eliminating their associated operators as well. This means that for the operators of the following layers that start their bit, the input they receive from $b$ will be the carry of the operator of the previous bit in the previous layer. Next, because any nontrivial factor of an $n$-bit integer is strictly smaller than $N$, the most significant bit of the exposed $p$-register can be set to $0$, eliminating its associated operator. Finally, since the hard instances considered in this work are odd, the least significant bit of both $p$ and $q$ has to be $1$. This implies in the case of $q$ that its first diagonal can be omitted, since it always passes $1$. The resulting logic circuit is shown in Fig.~\ref{fig: Logical Circuit} c.

\subsection{Tensor Network equation}\label{ssec: tensorization}
The next step is to tensorize the logical circuit, as explained in~\cite{melocoton}. This process involves transforming the operators into their tensor representations by replacing each operator with a tensor which has as many indices as the total number of inputs and outputs of the operator. The transformation for a generic tensor is shown in Fig.~\ref{fig: Tensorization} and the mathematical representations are defined in Appendix~\ref{appendix: tensor definitions}.

\begin{figure}
    \centering
    \includegraphics[width=\linewidth]{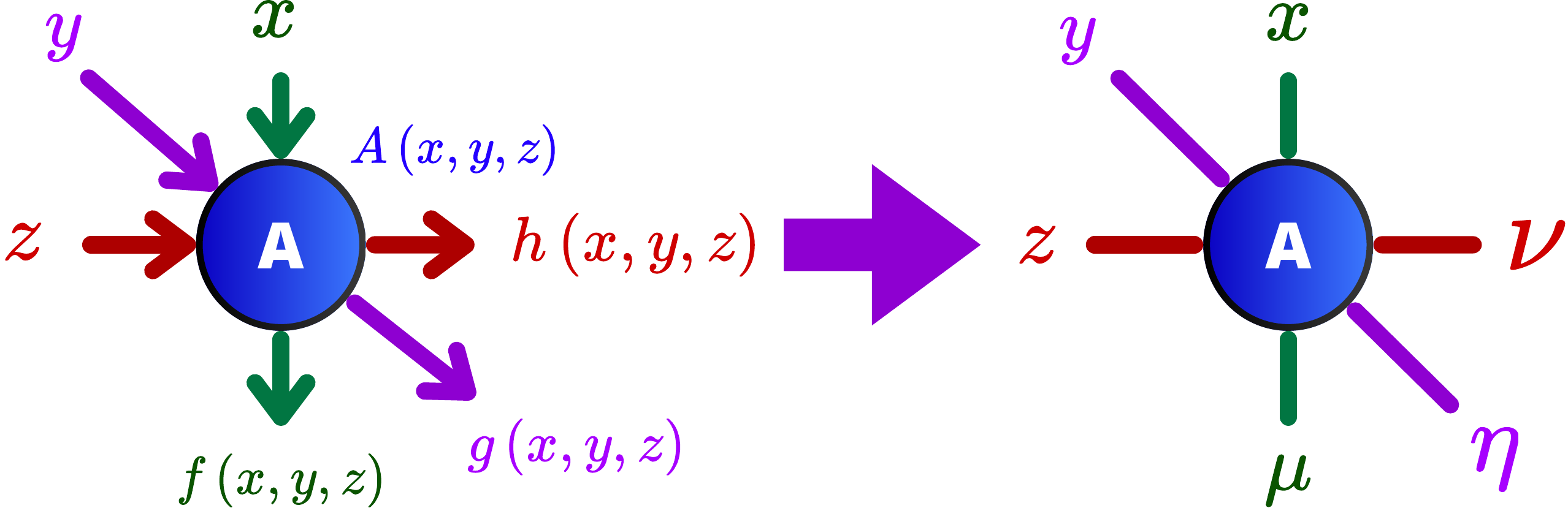}
    \caption{Tensorization of a six indexes operator with 3 inputs and 3 outputs that has an internal transformation $A_{x,y,z,\mu,\nu,\eta}=A\left(x,y,z\right)\delta_{\mu, f(x,y,z)}\delta_{\eta, g(x,y,z)}\delta_{\nu, h(x,y,z)}$. This means the resulting tensor has only non-zero elements $A_{x,y,z,\mu,\nu,\eta}=A\left(x,y,z\right)$ when $\mu=f(x,y,z)$, $\eta, g(x,y,z)$ and $\nu=h(x,y,z)$.}
    \label{fig: Tensorization}
\end{figure}

Tensorization considers at the same time all possible inputs $p$ and $q$, allowing the algorithm to determine the admissible values of the exposed register $p$. Let $p\in\{0,1\}^{n-1}$ be the exposed register, let $q\in\{0,1\}^{m}$ be the reduced auxiliary register with $m=\lceil n/2\rceil$, and let $y\in\{0,1\}^{W}$ be the output register, with $W$ large enough to represent the full product. In this formulation, the tensor network obtained after tensorizing the circuit represents a tensor $T(p,q,y)$.

\begin{theorem}[Circuit tensorization]
Let $C_{\times}$ be a deterministic Boolean circuit computing the binary product $y=pq$ on an output register of width $W$. Replace each gate by its indicator tensor and contract all internal wire indices. Then the resulting tensor satisfies
\begin{equation}
    T(p,q,y)=\mathbf{1}[y=pq].
\end{equation}
\end{theorem}
\textit{Proof.} For fixed external values $p$, $q$ and $y$, determinism implies that there is at most one assignment of the internal wires compatible with the circuit computation. If $y$ equals the product $pq$, that internal assignment satisfies every gate tensor and the contraction contributes exactly $1$. If $y$ is not the product $pq$, then every assignment violates at least one gate constraint, so the contraction is $0$. \hfill$\square$

The target integer $N$ is embedded into the same $W$-bit register by padding its binary representation with leading zeros. Projecting the output of the network onto this padded target produces the relation
\begin{equation}
    M_N(p,q)=\sum_{y} T(p,q,y)\prod_{i=0}^{W-1}\delta_{y_i,N_i}
    =\mathbf{1}[pq=N].
\end{equation}
Therefore, the projection enforces exact integer multiplication rather than equality modulo a power of two.

To encode the intended domain restrictions, let $B(p,q)$ be an admissibility selector. In the basic construction, $B(p,q)$ excludes trivial values and keeps the parity and register-length restrictions used in the odd semiprime setting. Additional selectors may be included if desired, but the reduced $q$-register is not itself an ordering condition. The exposed marginal is then
\begin{equation}
    S_N(p)=\sum_{q\in\{0,1\}^{m}} M_N(p,q)B(p,q).
\end{equation}
Its support is characterized by
\begin{equation}
    \operatorname{supp}S_N=\left\{p:\exists q<2^m \text{ such that } pq=N \text{ and } B(p,q)=1\right\}.
\end{equation}

\begin{theorem}[Reduced $q$-register and support]
Let $N$ be an $n$-bit composite integer and let $m=\lceil n/2\rceil$. Suppose the network enforces exact multiplication $pq=N$ and allows nontrivial auxiliary values $q$. Then $S_N$ has nonempty support. Moreover, every $p\in\operatorname{supp}S_N$ is a nontrivial divisor of $N$, provided the admissibility constraints exclude $p=1$, $p=N$, $q=0$ and $q=1$.
\end{theorem}
\textit{Proof.} Since $N$ is composite, it has a nontrivial divisor $d$ with $1<d\leq\sqrt{N}$. Because $\sqrt{N}<2^{n/2}\leq 2^m$, the value $d$ fits in the reduced $q$-register. Hence the assignment $q=d$ and $p=N/d$ is representable, so the support is nonempty. Conversely, if $p\in\operatorname{supp}S_N$, then there exists an admissible $q$ such that $pq=N$, and therefore $p$ divides $N$. The nontriviality assumptions remove the cases $p=1$ and $p=N$. \hfill$\square$

For a semiprime $N=ab$ with $1<a\leq b$, the reduced $q$-register guarantees that the complementary factor $b$ belongs to $\operatorname{supp}S_N$ through the assignment $q=a$. If $b$ also fits in the reduced register, the reverse orientation may also survive. Both exposed values are valid nontrivial factors, and the construction does not require identifying the smaller factor in advance.

From the tensor network, obtaining a factor can be achieved using the iterative Half Partial Trace~\cite{melocoton}. The direct scalar rule is exact when the current support is a singleton. If $\operatorname{supp}S_N=\{p^{\star}\}$, then
\begin{equation}
    S_N(p)=\delta_{p,p^{\star}},
\end{equation}
and for each bit $k$ the two-component marginal
\begin{equation}
    V_b^{(k)}=\sum_{p}\delta_{b,p_k}S_N(p),\qquad b\in\{0,1\},
\end{equation}
satisfies $V_b^{(k)}=\delta_{b,p_k^{\star}}$. In that case
\begin{equation}
    \Omega_k=V_1^{(k)}-V_0^{(k)}\in\{-1,1\},
\end{equation}
so the sign determines the bit exactly. To avoid any ambiguity at the origin, define
\begin{equation}
    H_{+}(t)=
    \begin{cases}
        1,& t>0,\\
        0,& t<0.
    \end{cases}
\end{equation}
Under the singleton-support hypothesis, $\Omega_k$ is never zero and the exact bit is recovered as $p_k^{\star}=H_{+}(\Omega_k)$.

If several admissible factors remain, then $V_b^{(k)}$ is a counting marginal over the current support rather than the bit value of a unique factor. In that multi-support case, independent decisions based only on the sign of $\Omega_k$ may mix bits coming from different factors. A correct reconstruction is instead obtained by adaptive projection. Starting from $S_0(p)=S_N(p)$, suppose that after $r$ steps one has a projected tensor $S_r(p)$ with nonempty support. Choose an unfixed bit $k$ and compute
\begin{equation}
    V_b^{(r,k)}=\sum_{p}\delta_{b,p_k}S_r(p),\qquad b\in\{0,1\}.
\end{equation}
Select any branch $b$ such that $V_b^{(r,k)}>0$, and project
\begin{equation}
    S_{r+1}(p)=S_r(p)\delta_{p_k,b}.
\end{equation}
Because the chosen branch has positive marginal, the new support remains nonempty. After all bits are fixed, the resulting bit string $p^{\star}$ belongs to the original support, and therefore it is a valid nontrivial factor.

\begin{figure*}
    \centering
    \includegraphics[width=\linewidth]{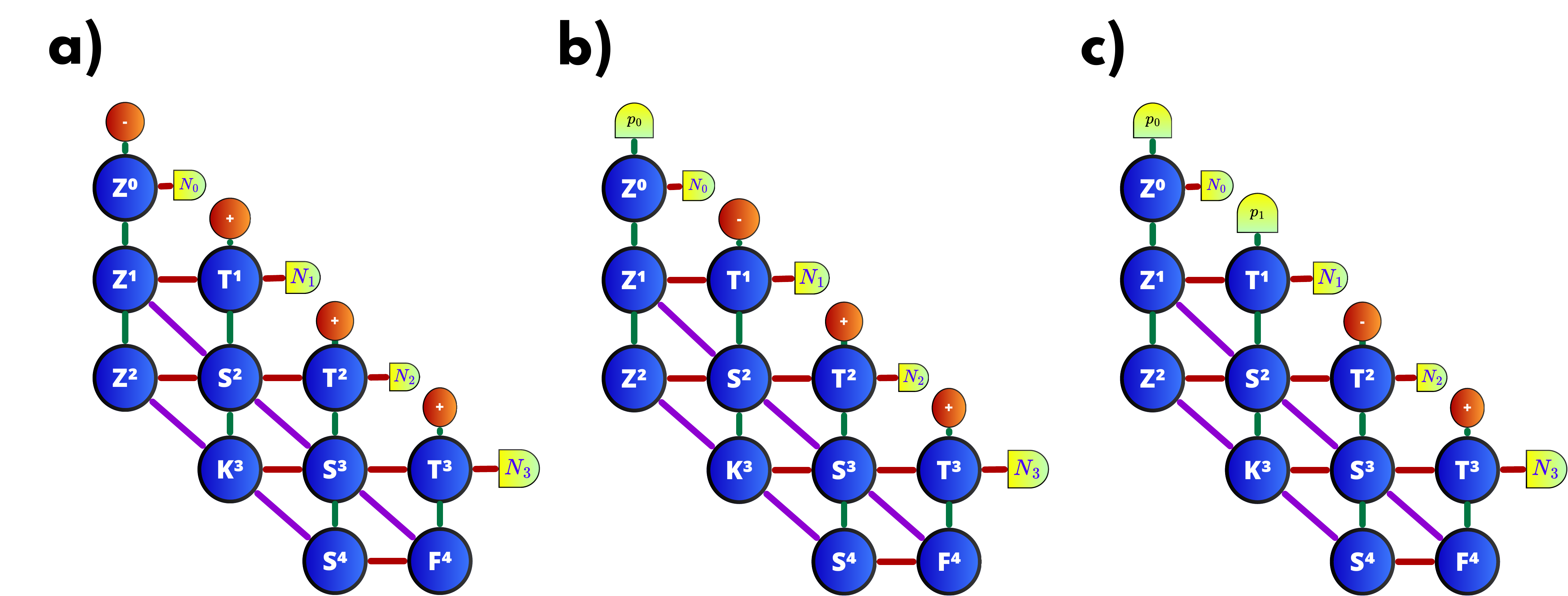}
    \caption{Tensor network of the prime factorization for five bits number $N$ which determines the a) first bit of $p$, b) second bit of $p$, c) third bit of $p$.}
    \label{fig: Tensor network}
\end{figure*}

The defined tensor network is illustrated in Fig.~\ref{fig: Tensor network}, where the previously selected values can be imposed when the support is reduced iteratively. It can be solved from the least significant bit to the most significant bit, as shown in Fig.~\ref{fig: Tensor network}, or from the most significant to the least. This will depend on the contraction scheme and the compression process used. In the exact singleton-support case this additional projection is not necessary, whereas in the multi-support case it is the mechanism that selects one admissible factor.

Once a candidate $p^{\star}$ is reconstructed, its validity can be checked by ordinary integer division. If $1<p^{\star}<N$ and $p^{\star}$ divides $N$, then $q^{\star}=N/p^{\star}$ gives the complementary factor. Under a semiprime promise, this completes the factorization. For a general composite integer, the same procedure can then be applied recursively to the factors obtained.

It is interesting to note that certain constraints can be removed from the equation by adding more complexity to the tensor network. The odd condition can be removed by introducing the first diagonal, allowing the value $0$ for the least significant bit of $q$. An optional selector $q\leq p$ can also be introduced to canonicalize the orientation while still allowing the square case $q=p=\sqrt{N}$. More general composite numbers can be treated by using the enforcement of previous results and by recursively applying the same nontrivial-divisor search to the factors obtained.

\subsection{Contraction scheme and computational complexity analysis}\label{ssec: complexity}

With the tensor network equation defined, the results of the equation can be computed with a contraction scheme, the order chosen for the contractions of the pairs of tensors of the network. In this section, the computational complexity of computing the matrix-matrix multiplication
\begin{equation}
    C = \sum_{i,j,k}^{I,J,K} A_{ij}B_{jk} \hat{e}_{i}\hat{e}_{k}
\end{equation}
is assumed to be $O(IJK)$.

\begin{figure*}
    \centering
    \includegraphics[width=\linewidth]{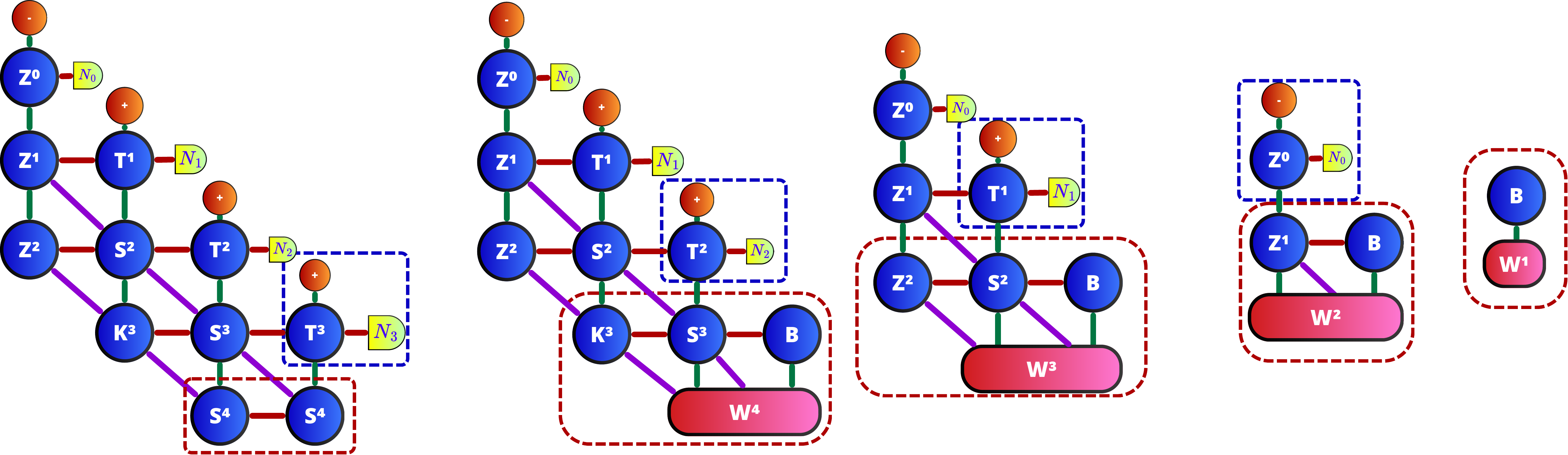}
    \caption{Exact contraction scheme from bottom-up. Last row is contracted, and the resulting tensor is contracted with the the row above it from right to left.}
    \label{fig: Down to up}
\end{figure*}

The tensor network contraction scheme can be performed mainly in two ways: from bottom to top and from left to right. In the bottom-up case, each tensor of the last row is contracted from right to left, until the whole row is contracted. Since in the last step, a tensor of $n-1$ indices of dimensions 
\begin{equation*}
    (\underbrace{2, \dots, 2}_{n/2-1}, \underbrace{3, \dots, 3}_{n/2-1}, 4)
\end{equation*}
is contracted with one of three indices of dimensions $(4,3,2)$, through its dimension $4$ index, the computational complexity of contracting this row is \mbox{$O(4(2\cdot3)^{n/2})=O(6^{n/2})$}. After that, this tensor is contracted with each tensor of the next row, from right to left, until the entire row is fully contracted. In this case, each contraction has a cost of $O(6^{n/2})$, and for the whole row, $O(n6^{n/2})$. Applying this contraction process to all rows results in an overall complexity of $O(n^26^{n/2})$. Moreover, since this procedure must be repeated $O(n)$ times in the iterative process to determine each bit of $p$, the total cost is $O(n^3 6^{n/2})$. In terms of $N$, the computational complexity is $O\left((\log_2 N)^3 6^{\log_2 \sqrt{N}}\right)$ or $O\left((\ln N)^3 \exp\left(\left(\log_2 \sqrt{6}\right) \ln N\right)\right)$. This complexity is worse than the brute force complexity. The contraction scheme is shown in Fig.~\ref{fig: Down to up}.

\begin{figure*}
    \centering
    \includegraphics[width=\linewidth]{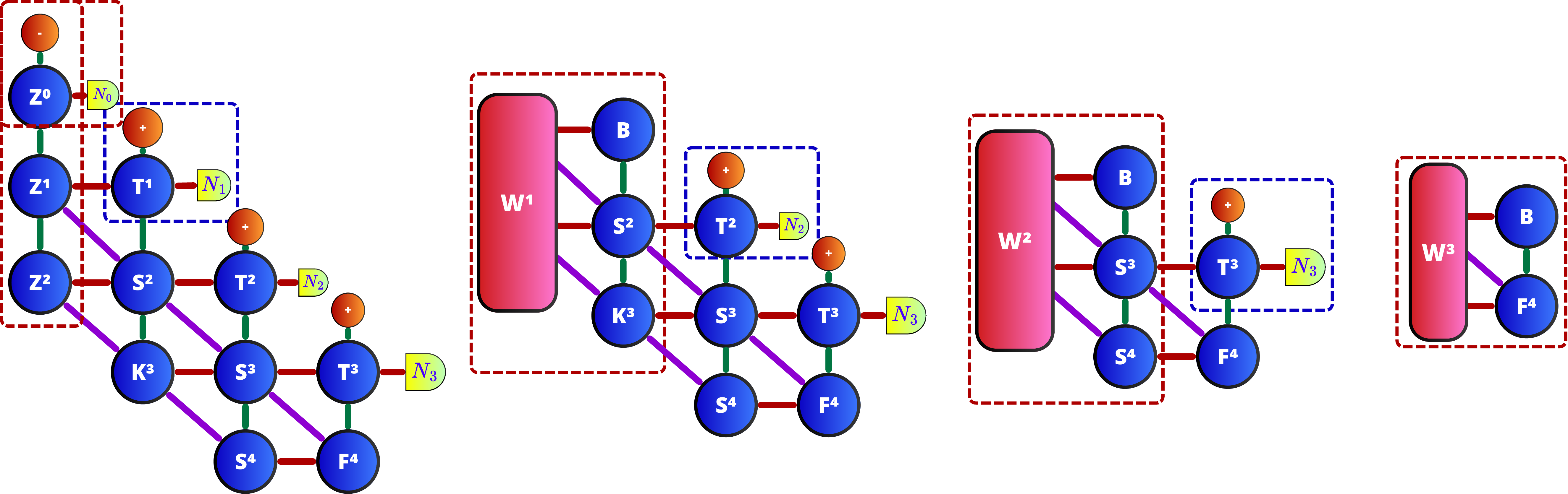}
    \caption{Exact contraction scheme from left-right. First column is contracted, and the resulting tensor is contracted with the the row above it from top to bottom.}
    \label{fig: left-right}
\end{figure*}

In the left-to-right case, it is contracted following the same style. First, the initial column is contracting from top to bottom and then with the following columns from top to bottom. Since the tensors have dimension $O(4^{n/2})$, the cost of contracting each column is $O(n4^{n/2})$. Repeating this for $n$ columns results in a complexity of $O(n^2 4^{n/2})$ per iteration, which to obtain the $n$ bits of $p$ results in a total cost of $O(n^3 4^{n/2})$. In terms of $N$, the computational complexity is $O\left((\log_2 N)^3 4^{\log_2 \sqrt{N}}\right)$ or $O\left((\ln N)^3 N\right)$. This complexity is also worse than the brute force complexity. The contraction scheme is shown in Fig.~\ref{fig: left-right}.

\paragraph{Sparse-aware left-to-right contraction.}
The previous left-to-right estimate treats the intermediate tensors as dense arrays over their ambient frontier space. However, the local tensors listed in Appendix~\ref{appendix: tensor definitions} are indicator tensors of constant-arity logical relations, so the same tensor network and the same left-to-right contraction order can be evaluated exactly with sparse storage. In that representation, each intermediate tensor is stored as a table of its non-zero entries: the key is the current assignment of the frontier indices and the value is the accumulated multiplicity. Contracting one more column then amounts to sparse joins on the shared indices, followed by the summation or projection of the eliminated ones. The tensor network and the contraction order are unchanged; only the representation of the intermediate tensors changes from dense arrays to sparse tables of non-zero frontier assignments.

Let $\mathcal{R}_r$ denote the intermediate sparse tensor obtained after contracting the first $r$ columns from left to right, that is, after imposing the first $r$ least significant output constraints of the exact product relation $pq=N$. The following argument therefore refers to the low-to-high version of the left-to-right order. Since the relevant hard instances are odd, and an even $N$ is handled by first extracting the factor $2$, one has $p_0=q_0=1$ in the odd case considered here. Every non-zero entry of $\mathcal{R}_r$ must then satisfy
\begin{equation}
    pq\equiv N \pmod{2^r}.
\end{equation}
For every admissible odd $q$-prefix of length $\min(r,m)$, the corresponding residue class is invertible modulo $2^r$, so there is a unique compatible $p$-prefix satisfying
\begin{equation}
    p\equiv Nq^{-1} \pmod{2^r}.
\end{equation}
Moreover, in the optimized deterministic carry encoding of Appendix~\ref{appendix: tensor definitions}, once the active prefixes of $p$ and $q$ and the first $r$ output bits are fixed, the compatible frontier carry states are determined by the local update relations. Consequently,
\begin{equation}
    \operatorname{nnz}(\mathcal{R}_r)=O\left(2^{\min(r,m)}\right),
\end{equation}
with at most a constant-factor modification if a different but equivalent carry encoding is used. In particular,
\begin{equation}
    \max_r \operatorname{nnz}(\mathcal{R}_r)=O(2^m).
\end{equation}
This characterizes the non-zero support of the same exact contraction and does not reformulate the method as an external enumeration of divisors.

Each column contains $O(m)$ local tensors of constant arity and constant index dimension. Therefore, using sparse hash-join contractions, the cost of absorbing one more column is
\begin{equation}
    O\left(m\,\operatorname{nnz}(\mathcal{R}_r)\right).
\end{equation}
Summing over the $n$ columns gives the cost of one sparse left-to-right contraction,
\begin{equation}
    T_{\mathrm{sparse,LTR}}
    =
    \sum_{r=1}^{n} O\left(m\,\operatorname{nnz}(\mathcal{R}_r)\right)
    =
    O(nm2^m)
    =
    O\left(n^2 2^{\lceil n/2\rceil}\right).
\end{equation}
In terms of $N$, this can be written as $O((\log N)^2\sqrt{N})$ up to constant factors and bit-complexity details. Hence sparse storage changes the exponential dependence of the left-to-right bound from the dense ambient frontier scale $4^m$ to the exact support scale $2^m$, improving the exponential part of the dense estimate from $O(n^2 4^{n/2})$ per contraction to $O(nm2^m)$.

The corresponding streaming memory requirement is $O(2^m)$ sparse entries when only one forward message is stored at a time. If one also counts the frontier keys explicitly, this is $O(m2^m)$ stored bits up to constant factors; if forward and backward messages are cached simultaneously for reuse, the memory increases to $O(n2^m)$ sparse entries up to the same key-encoding factors.

For the iterative Half Partial Trace there are two sparse regimes. If each bit marginal is recomputed from scratch, the conservative total cost becomes
\begin{equation}
    O\left(n^3 2^{\lceil n/2\rceil}\right),
\end{equation}
because there are $O(n)$ contractions and each one costs $O(nm2^m)=O(n^2 2^m)$. In contrast, the exact iterative projection procedure described above imposes the already selected bits of $p$ as delta constraints. After $k$ low-order bits of $p$ have been fixed and the support remains nonempty, the same congruence relation restricts the compatible low-order prefixes of $q$ to $O(2^{m-k})$. Therefore the total cost of the exact iterative sparse contraction is bounded by
\begin{equation}
    \sum_{k=0}^{m} O(nm2^{m-k}) + O((n-m)nm)
    =
    O(nm2^m)
    =
    O\left(n^2 2^{\lceil n/2\rceil}\right).
\end{equation}
This estimate applies to the exact iterative projection procedure, where each selected bit value is imposed as a delta constraint and the remaining support is kept nonempty. The sparse bound remains exponential in $n$ and should not be interpreted as a polynomial-time factoring algorithm. As in the dense baseline, the one-time explicit creation of all tensor instances would still require $O(n^2)$ work.
Under this exact iterative sparse procedure, the total left-to-right cost therefore decreases from the dense baseline $O(n^3 4^{n/2})$ to $O\left(n^2 2^{\lceil n/2\rceil}\right)$.

In all of these scenarios, the time complexity and the required memory scale exponentially in the number of bits of $N$. In an attempt to overcome this limitation, the possibility of computing this equation approximately is evaluated. The complexity is due to the high size of the tensors in the intermediate computations. The initial and final idealized tensors may have simple product-like structure, because both the input of the problem $N$ and any singleton output $p$ are products and can therefore be represented by a matrix product state (MPS), or tensor train (TT), of bond dimension $1$. However, this does not imply that the intermediate TT ranks remain small: they can grow substantially during the contraction. Therefore, the observed compressibility of the intermediates is an empirical property of the tested instances and contraction scheme, not a general consequence of Rolle's theorem. Dense complexity estimates should thus be interpreted as baselines, while sparse-aware contractions and TT/MPS compression are practical accelerations whose performance depends on the contraction order, intermediate ranks, truncation thresholds and numerical stability.

This point is particularly important for the scalar quantities used in the bit-selection stage. In the exact singleton-support case one has $\Omega_k\in\{-1,1\}$, so for an approximation $\widetilde{\Omega}_k$ the sign is guaranteed to be preserved whenever $|\widetilde{\Omega}_k-\Omega_k|<1/2$. In the multi-support case, however, the relevant objects are the positive marginals $V_b^{(r,k)}$, and a reliable approximate computation must preserve the distinction between a genuinely empty branch and a positive branch. In all cases, any reconstructed candidate factor must be validated afterward by exact division.

\section{Performance evaluation and results}\label{sec: performance}
In this section, the correctness of the presented equation and the performance of its calculation are evaluated. The exact case is implemented by means of bottom-up contraction, to better evaluate the computational time scaling in Sec.~\ref{ssec: exact}. In Sec.~\ref{ssec: approx}, the approximated case is studied, and the minimum bond dimension needed to obtain the correct result is determined depending on the problem size. Finally, it is studied whether there is a dependence of the number of erroneous bits in the solution obtained on the ratio between the bond size used and that required in the failed cases.

For all these tests, ten numbers are used to factorize for each value of $n$ (if available), so that the prime numbers used for the construction of $N$ are not excessively different. This is to take the worst case, since if one is much smaller than another, the problem is simpler to solve. The runs have been performed on a regular laptop CPU Intel(R) Core(TM) i7-14700HX 2.10 GHz, with 32GB RAM and 20 cores, using the Python Tensorkrowch~\cite{TK} library.

\subsection{Exact case}\label{ssec: exact}
The first step is to compute the result of the equation with the presented contraction schemes. The contraction is performed using a bottom-up scheme. In addition, the time $t_{total}$ is evaluated, defined as the time needed to obtain the solution from the moment the number $N$ is introduced into the function to the moment the result is retrieved and the contraction time $t_{contr}$. No early stopping or bit prefixing at $0$ is implemented to avoid biasing the evaluation.

\begin{figure}
    \centering
    \includegraphics[width=\linewidth]{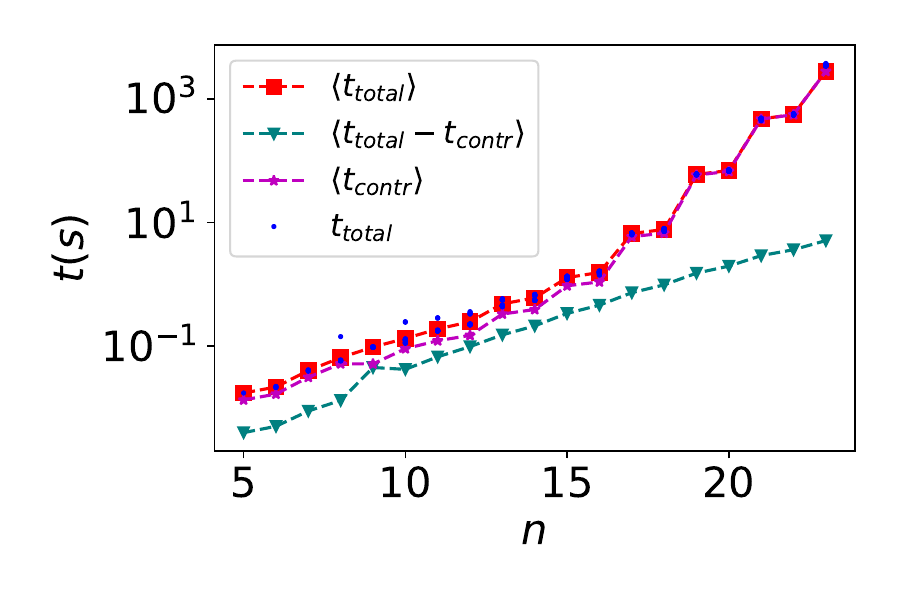}
    \caption{Execution time vs number $n$ of bits of $N$. In blue points the time for each instance and in red squares the average total time. The purple star are the contraction time and the green triangles are the rest of time.}
    \label{fig: time vs n}
\end{figure}

As expected, in every attempt, the equation provides a valid exposed factor, and in the tested semiprime convention it matches the target register $p$. Fig.~\ref{fig: time vs n} shows the total time of each execution, the average for each problem size, the average contraction time, and the average of the time which is not due to the contraction. As previously predicted, execution time scales exponentially with the $n$ number of bits of $N$, and its main contribution is the contraction time. However, the exponential scaling of the time not related to the contraction may indicate that the implementation of the creation of the tensors and the processing of the data could be improved. These implementation problems may explain the slightly ascending curvature of the contraction time. It is also shown that the times do not have a relevant dispersion, probably caused by its execution on a regular laptop.

\subsection{Approximated case}\label{ssec: approx}
\begin{figure}[h]
    \centering
    \begin{subfigure}[b]{0.48\linewidth}
        \centering
        \includegraphics[width=\linewidth]{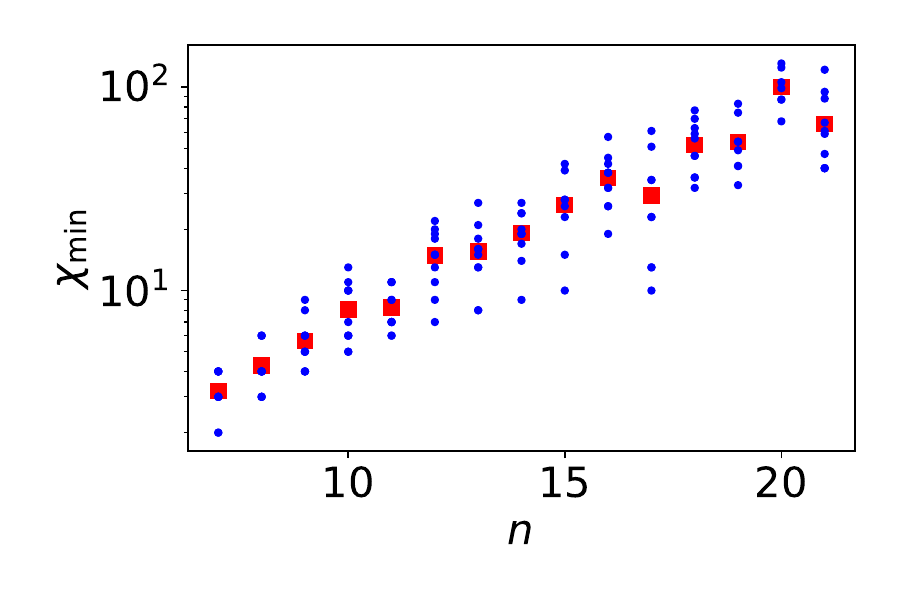}
        \caption{Minimal bond dimension $\chi_{\min}$ vs number $n$ of bits of $N$ for bottom-up scheme. In blue points the dimension for each instance and in red squares the average dimension.}
        \label{fig: bond down}
    \end{subfigure}
    \hfill
    \begin{subfigure}[b]{0.48\linewidth}
        \centering
        \includegraphics[width=\linewidth]{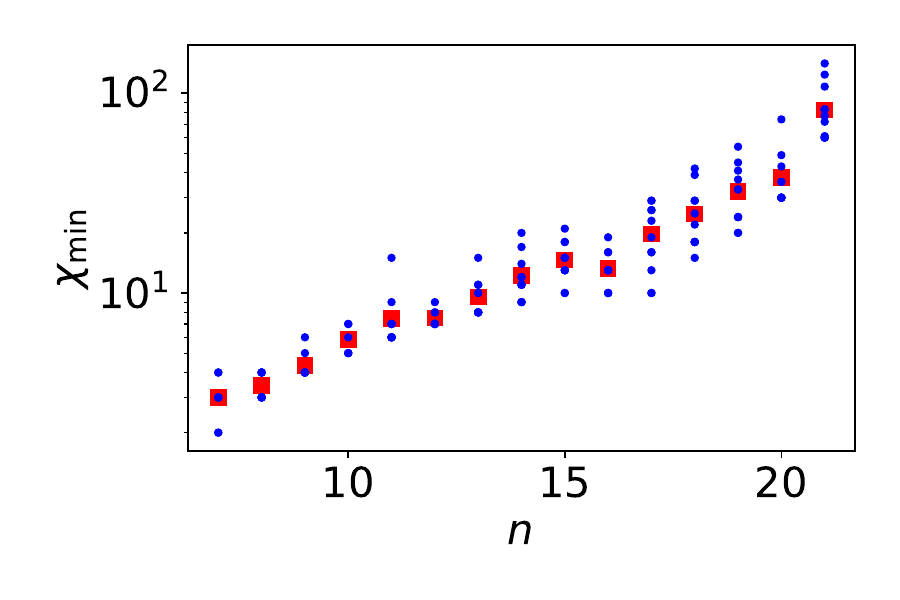}
        \caption{Minimal bond dimension $\chi_{\min}$ vs number $n$ of bits of $N$ for left-right scheme. In blue points the dimension for each instance and in red squares the average dimension.}
        \label{fig: bond left}
    \end{subfigure}
    \caption{Minimal bond dimension $\chi_{\min}$ vs number $n$ of bits of $N$ for different contraction schemes.}
    \label{fig:bond_dimensions}
\end{figure}

Now, the focus is to compute it more efficiently using an approximate train tensor. Figs.~\ref{fig:bond_dimensions} show the minimum bond dimension needed to obtain the correct solution for each instance, and the average for each value of $n$. Fig.~\ref{fig: bond down} shows the bottom-up contraction, and Fig.~\ref{fig: bond left} shows the left-right contraction.

\begin{figure}[h]
    \centering
    \begin{subfigure}[b]{0.48\linewidth}
        \centering
        \includegraphics[width=\linewidth]{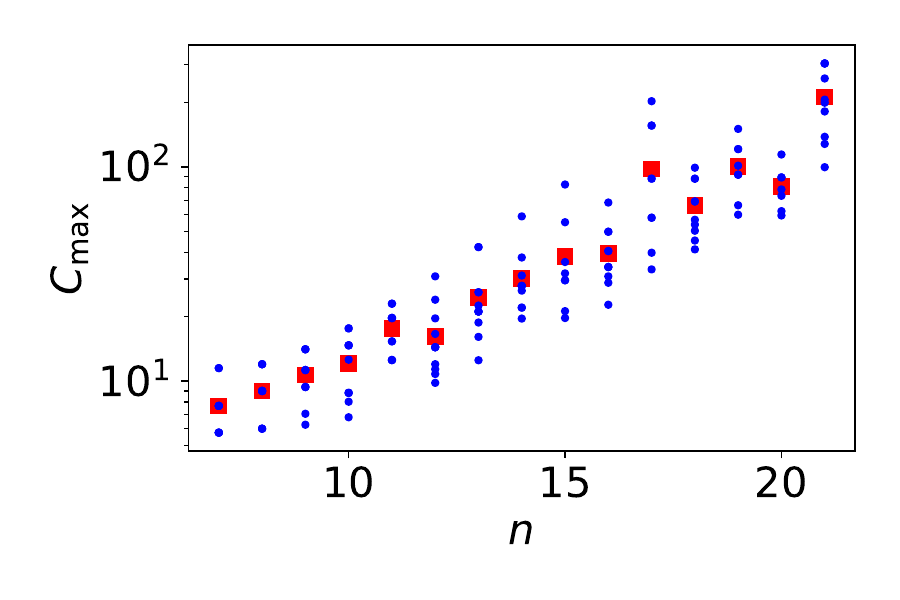}
        \caption{Maximum possible compression vs number $n$ of bits of $N$ for bottom-up scheme.}
        \label{fig: compression down}
    \end{subfigure}
    \hfill
    \begin{subfigure}[b]{0.48\linewidth}
        \centering
        \includegraphics[width=\linewidth]{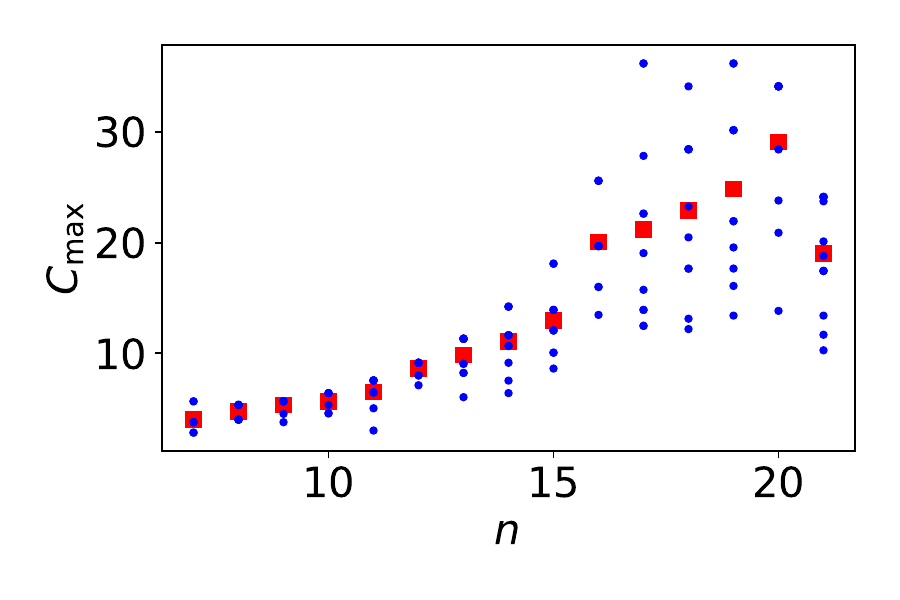}
        \caption{Maximum possible compression vs number $n$ of bits of $N$ for left-right scheme.}
        \label{fig: compression left}
    \end{subfigure}
    \caption{Maximum possible compression vs number $n$ of bits of $N$ for different contraction schemes.}
    \label{fig:compression_rates}
\end{figure}

It can be observed that the minimum bond dimension grows exponentially with the number of bits of the number to factorize. This implies that the algorithm continues to scale exponentially with the size of the problem. However, the scaling is much smaller than in the exact case since it arrives at bond dimensions much smaller than what would be expected in the exact case. The minimum bond dimension in the exact uncompressed case is
\begin{equation}
    r_{\max} = \sqrt{6^{n/2}},\qquad  r_{\max} = \sqrt{4^{n/2}},
\end{equation}
for the bottom-up scheme and the left-right scheme, respectively, so the compression factor is defined as
\begin{equation}
    C_{\max} = \frac{r_{\max}}{\chi_{\min}}.
\end{equation}
As can be seen in Figs.~\ref{fig:compression_rates}, the scaling of the compression factor is also exponential with respect to $n$. This indicates that there is not as much entanglement in the system as was initially considered and that there may be a more efficient tensor network to solve the problem.

\begin{figure}[h]
    \centering
    \begin{subfigure}[b]{0.48\linewidth}
        \centering
        \includegraphics[width=\linewidth]{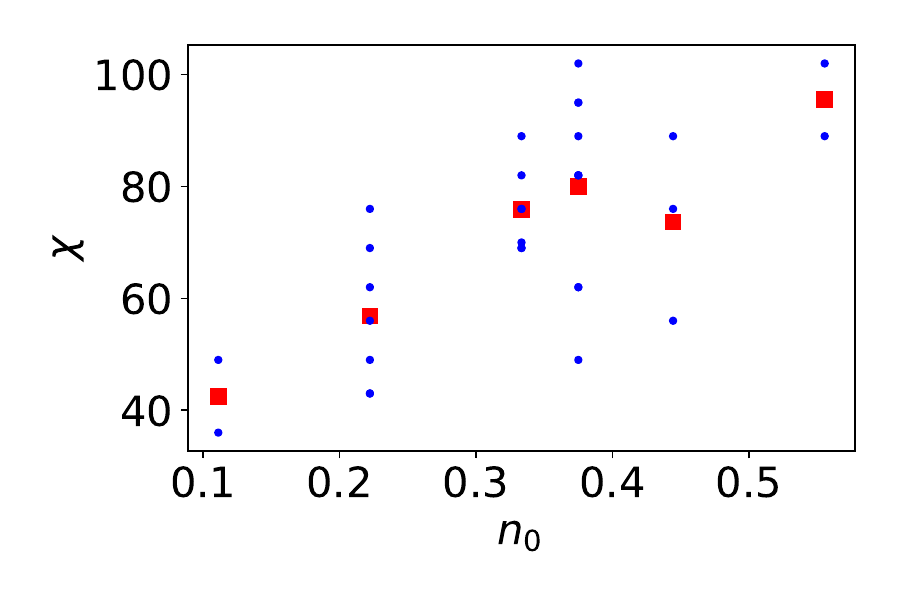}
        \caption{Minimal bond dimension vs $n_0$ for $n=16$ with $p=179$ for bottom-up scheme.}
        \label{fig: dimension zeros down}
    \end{subfigure}
    \hfill
    \begin{subfigure}[b]{0.48\linewidth}
        \centering
        \includegraphics[width=\linewidth]{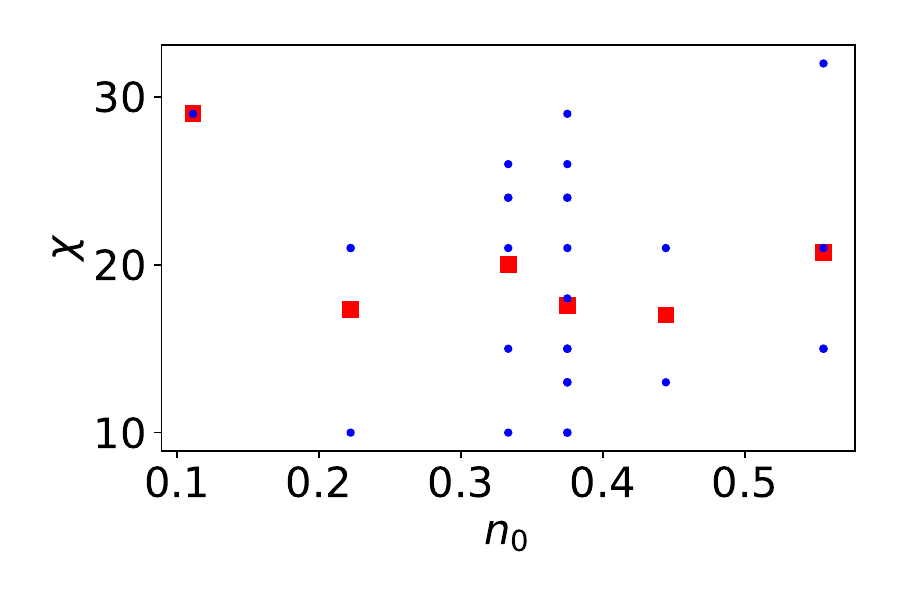}
        \caption{Minimal bond dimension vs $n_0$ for $n=16$ with $p=179$ for left-right scheme.}
        \label{fig: dimension zeros left}
    \end{subfigure}
    \caption{Minimal bond dimension vs $n_0$ for different contraction schemes.}
    \label{fig:dimension_zeros}
\end{figure}

An interesting point of study is to check if there is a dependence between the required bond dimension and the number of zeros in one of the prime factors. The proportion of zeros is defined as
\begin{equation}
    n_0 = \frac{n(p,0)}{n(p,0)+n(p,1)},
\end{equation}
being $n(x,y)$ the number of bits with value $y$ in binary form of $x$. The results of Figs.~\ref{fig:dimension_zeros} show a relation cannot be directly deduced in the left-right scheme, but it seems that there may be a dependency in the bottom-up. However, given the sample, no relation can be concluded.

Finally, it is checked if there is a possibility that, although the bond dimension used is lower than the minimum necessary, the result obtained allows us to obtain part of the information of the correct solution. The bit-flip error is defined as
\begin{equation}
    B = \sum_{i=0}^{n-1} \frac{(1-x_i)p_i + x_i(1-p_i)}{n}.
\end{equation}

\begin{figure}
    \centering
    \begin{subfigure}[b]{0.48\linewidth}
        \centering
        \includegraphics[width=\linewidth]{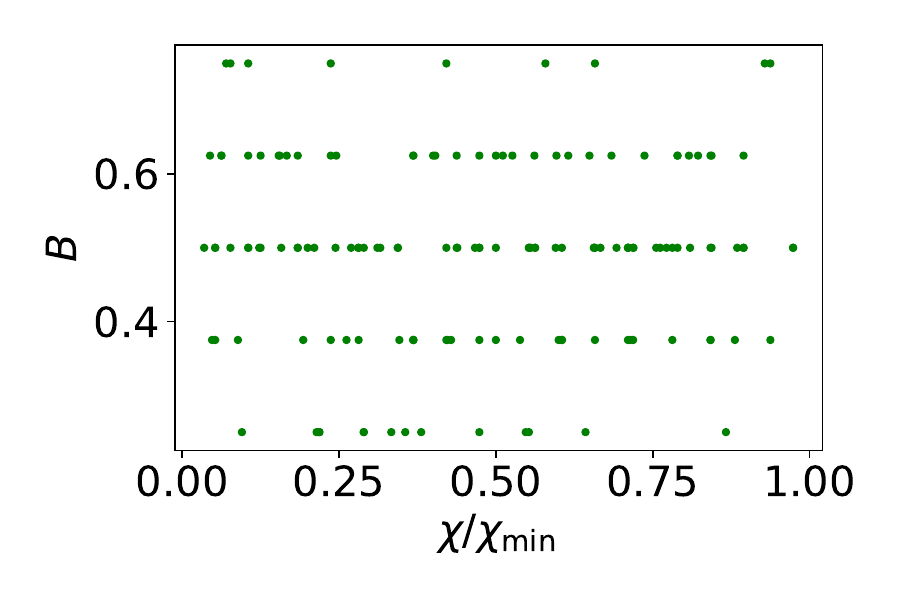}
        \caption{Bit-flip rate $B$ vs the proportion $\frac{\chi}{\chi_{\min}}$ for bottom-up scheme.}
        \label{fig: error down}
    \end{subfigure}
    \hfill
    \begin{subfigure}[b]{0.48\linewidth}
        \centering
        \includegraphics[width=\linewidth]{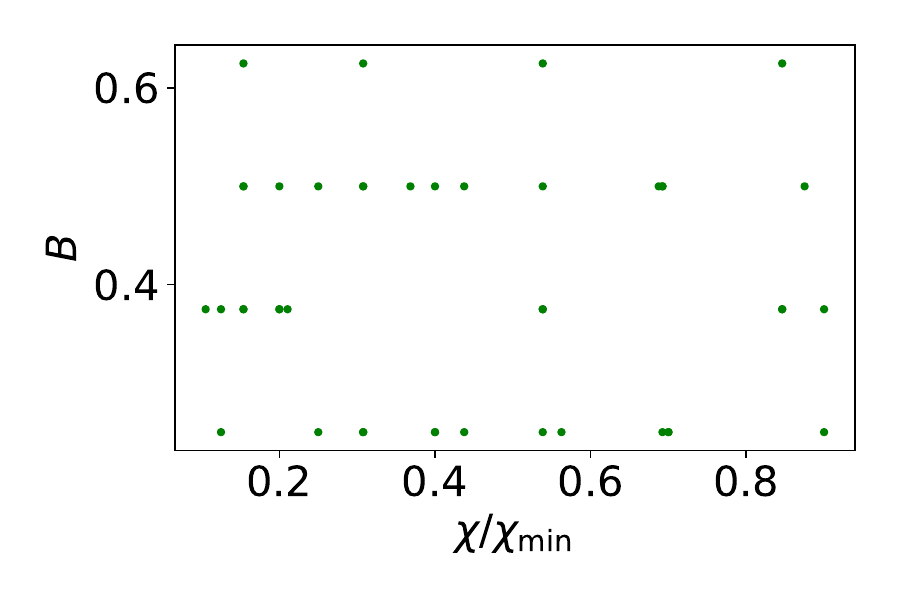}
        \caption{Bit-flip rate $B$ vs the proportion $\frac{\chi}{\chi_{\min}}$ for left-right scheme.}
        \label{fig: error left}
    \end{subfigure}
    \caption{Bit-flip rate $B$ vs the proportion $\frac{\chi}{\chi_{\min}}$ for different contraction schemes.}
    \label{fig:error_rates}
\end{figure}

As can be seen in Figs.~\ref{fig:error_rates} , there does not appear to be any relationship between the level of error and the increase in bond dimension.

\section{Discussions}\label{sec: discussions}

This work has presented an exact and explicit equation that allows recovering a nontrivial factor of an odd semiprime number, and more generally encodes the search for nontrivial divisors of a composite integer. However, its computation requires an exponential amount of time and memory, making its use worse than the brute force search for these factors. This equation has been computed with classical computational resources and its correctness and exponential scaling have been verified. In the approximated computation, tensor train compression has been studied to improve its performance, but the results remain exponential. Nevertheless, the results obtained with high compression indicate that there could exist a more efficient tensor network equation, or a more efficient contraction strategy, for the same relation. This is consistent with the high sparsity of all the tensors, even the intermediate ones. The absence of a clear relation between the number of zeroes and the minimal bond dimension is certainly unexpected, because more equal digits may require simpler rules to implement. In addition, the absence of a clear relation between the bit-flip error and the bond dimension below the minimum required dimension is an interesting result. This may indicate that the system must either have the necessary bond dimension or it does not allow any useful information to be extracted from the correct solution. Both observations may be consequences of the tensor train representation, because its natural expressibility remains in one-dimensional relations, and the bit information does not correspond to a one-dimensional phenomenon. This opens the possibility of studying other tensor network representations for the compressed contraction.

\section{Conclusions}\label{sec: conclusions}

Even with these computational results, this work has provided an exact and explicit tensor-network equation for the computation of nontrivial factors of a target composite number. Its computation with tensor network contraction schemes has been studied, but there could be other more efficient ways to compute it. For example, if there existed a physical system that could efficiently compute this equation, it could be used to recover nontrivial factors of semiprime numbers and, under the semiprime promise, complete the factorization. This opens a new future research line in analog computing. Another possible way is to use mathematical properties of the equation to simplify it. For example, the equation is the contraction of several Kronecker delta tensors, and they can be decomposed into the contraction of other tensors. In this way, it could be applied to avoid certain computations. This opens up new research lines in pure mathematics and computational mathematics. Another possible research line from this work is the study of other tensor representations for compressed contraction, making use of the sparsity to improve the complexity and obtain a new tensor network that is easier to contract. Finally, this work provides a new perspective on how to approach this problem, and similar ones, with a simple tensor network construction.

\section*{Data and code availability}
All data and code required for this project can be accessed upon reasonable request by contacting the authors.

\section*{Acknowledgements}

\paragraph{Funding information}
The research has been funded by the Ministry of Science and Innovation and CDTI under ECOSISTEMAS DE INNOVACIÓN project ECO-20241017 (EIFEDE) and ICECyL (Junta de Castilla y León) under project CCTT5/23/BU/0002 (QUANTUMCRIP). This project has been funded by the Spanish Ministry of Science, Innovation and Universities under the project PID2023-149511OB-I00, and under the programme for mobility stays at foreign higher education and research institutions "José Castillejo Junior" with code CAS23/00340.

\begin{appendix}
\numberwithin{equation}{section}

\section{Tensor network definition}\label{appendix: tensor definitions}
There are the following types of tensors that make up the tensor network, to factorize an $N$ number of $n$ bits, allowing $p$ to have up to $n-1$ bits and $q$ to have up to $\lceil n/2 \rceil$ bits. In the exact formulation, the product is understood on an output register wide enough to contain the full integer product, with the bits above the $n$-bit expansion of $N$ projected to zero. The explicit tensors listed below correspond to the optimized network that remains after eliminating wires fixed by those boundary projections and by the parity constraints of the odd semiprime setting. Throughout this appendix, expressions such as $\mu=jl$ denote binary concatenation of index values rather than ordinary multiplication; similarly, dimension-$3$ indices are compound carry/output indices decoded by the formulas that accompany each tensor, and dimension-$4$ indices encode two binary degrees of freedom. Free horizontal indices correspond to bits of $p$, diagonal free indices correspond to bits of $q$, the output indices projected onto the target correspond to bits of $N$, and the remaining indices are contracted carry or partial-product indices.

\paragraph{Initial layer}

The first column is composed of three types of tensors, the initial $Z^0_{2\times 2\times 2}$ that connects to bit $p_0$ and bit $N_0$,
\begin{equation}
    \begin{gathered}
        \nu = \mu = j,\\
        Z^0_{j\mu\nu} = 1,
    \end{gathered}
\end{equation}
 the intermediate $Z^k_{2\times 2\times 2\times 2}$ for the intermediate rows,
\begin{equation}
    \begin{gathered}
        \nu = j,\\
        \mu = jl,\\
        Z^k_{jl\mu\nu} = 1,
    \end{gathered}
\end{equation}
and the last $Z^{\lceil n/2 \rceil}_{2\times 2\times 4}$ for the last row of the layer, carrying the most significant bit of $q$,
\begin{equation}
    \begin{gathered}
        \mu = jl,\\
        Z^{\lceil n/2 \rceil}_{j\mu l} = 1.
    \end{gathered}
\end{equation}

\paragraph{Presaturation layers}

The columns before saturation are composed of three types of tensors, the initial $T^k_{2\times 2\times 2\times 3}$ that connects to bit $p_k$ and bit $N_k$,
\begin{equation}
    \begin{gathered}
        \mu = i\oplus j,\\
        \nu = j\left(i+1\right),\\
        T^k_{ij\mu\nu} = 1,
    \end{gathered}
\end{equation}
the intermediate $S^k_{2\times 3\times 2\times 2\times 3 \times 2}$, for the intermediate rows,
\begin{equation}
    \begin{gathered}
        c =  \lceil j/2\rceil,\ y =  j - c,\\
        \mu = i\oplus c(l\oplus y),\\
        \nu = c\left(1+\left\lfloor \frac{i+l+y}{2}\right\rfloor\right),\\
        \eta = l,\\
        S^k_{ijl\mu\nu\eta} = 1,
    \end{gathered}
\end{equation}
and the last $K^k_{3\times 4\times 2\times 4}$, for the last row of the layer, carrying the most significant bit of $q$,
\begin{equation}
    \begin{gathered}
        c =  \lceil j/2\rceil,\ y =  j - c,\\
        y_i = \lfloor l/2\rfloor,\ y_l = l \ (\text{mod } 2),\\
        \mu = y_i\oplus c(y_l\oplus y),\\
        \eta = 2\left\lfloor \frac{y_i+c(y_l+y)}{2}\right\rfloor+y_l,\\
        K^k_{jl\mu\eta} = 1.
    \end{gathered}
\end{equation}

\paragraph{Saturation Layer}

The column of saturation is composed of three types of tensors, the initial $T^k_{2\times 2\times 2\times 3}$ that connects to bit $p_k$ and bit $N_k$ as before, the intermediate $S^k_{2\times 3\times 2\times 2\times 3 \times 2}$, for the intermediate rows as before, and the last $S^{n-1}_{3\times 4\times 4}$, for the last row of the layer, carrying the most significant bit of $q$ and informing the following layers if its bit has been in $1$.

Their non-zero elements are
\begin{equation}
    \begin{gathered}
        c =  \lceil j/2\rceil,\ y =  j - c,\\
        y_i = \lfloor l/2\rfloor,\  y_l = l \ (\text{mod } 2)\\
        \left\lceil \frac{y_i+c(y_l+y)}{2}\right\rceil= 0,\\
        \mu = 2y_l + y_i\oplus c(y_l\oplus y),\\
        S^{n-1}_{jl\mu} = 1.
    \end{gathered}
\end{equation}

\paragraph{Post-saturation layer}

The columns after saturation are composed of three types of tensors, the initial $T^k_{2\times 2\times 2\times 3}$ that connects to bit $p_k$ and bit $N_k$ as before, the intermediate $S^k_{2\times 3\times 2\times 2\times 3 \times 2}$, for the intermediate rows as before, and the last $R^{n-1}_{4\times 3\times 2\times 4}$, for the last row of the layer, carrying the most significant available bit of $q$ and informing the following layers if its bit has been in $1$.

Its non-zero elements are
\begin{equation}
    \begin{gathered}
        c =  \lceil j/2\rceil,\ y =  j - c,\\
        y_s = \lfloor i/2\rfloor,\ y_i = i \ (\text{mod } 2),\\
        c y_s (l+y) = 0, \left\lceil \frac{y_i+c(l+y)}{2}\right\rceil= 0,\\
        \mu = 2l + y_i\oplus c(l\oplus y),\\
        R^{n-1}_{ijl\mu} = 1.
    \end{gathered}
\end{equation}

\paragraph{Final layer}

The final column after is composed of two types of tensors, the initial $T^{n-2}_{2\times 2\times 2\times 3}$ that connects to bit $p_{n-2}$ and bit $N_{n-2}$ as before, and the last tensor $F^{n-1}_{4\times 3\times 2}$, for the last row of the layer, which imposes the final bit of $N$ to be $1$.

Its non-zero elements are
\begin{equation}
    \begin{gathered}
        c =  \lceil j/2\rceil,\ y =  j - c,\\
        y_s = \lfloor i/2\rfloor,\ y_i = i \ (\text{mod } 2),\\
        \text{if } y_s = 1\Rightarrow c(l+y) = 0,\\
        \left\lceil \frac{y_i+c(l+y)}{2}\right\rceil= 0,\\
        y_i\oplus c(l\oplus y) = 1,\\
        F^{n-1}_{ijl} = 1.
    \end{gathered}
\end{equation}
\end{appendix}



\bibliographystyle{unsrt}  
\bibliography{references}

\end{document}